# Forward Electricity Contract Price Projection: A Market Equilibrium Approach

Mateus A. Cavaliere, Sergio Granville, Gerson C. Oliveira, Mario V.F. Pereira

April 2019

*Abstract*--**This work presents a methodology for forward electricity contract price projection based on market equilibrium and social welfare optimization. In the methodology supply and demand for forward contracts are produced in such a way that each agent (generator/load/trader) optimizes a risk adjusted expected value of its revenue/cost. When uncertainties are represented by a discrete number of scenarios, a key result in the paper is that contract price corresponds to the dual variable of the equilibrium constraints in the linear programming problem associated to the optimization of total agents´ welfare. Besides computing an equilibrium contract price for a given year, the methodology can also be used to compute the evolution of the probability distribution associated to a contract price with a future delivery period; this an import issue in quantifying forward contract risks. Examples of the methodology application are presented and discussed.**

## I. INTRODUCTION

Electricity markets presents a great level risk due to spot prices uncertainty and volatility.

Volatility in electricity spot prices occurs in the short, medium and long-term and affects the cash flow of market agents who sell/buy energy in the wholesale market.

Short-term volatility is basically due to equipment (generators and transmission) failure, fuel availability, uncertainties in renewable energy production, fuel price and daily temperature. In medium and long-term periods, volatility in spot prices is due to uncertainty in hydrological conditions (in system with large reservoirs), demand growth, structural changes in the electricity generation sector and so on.

Figure 1 shows the hourly spot price on a randomly selected day in the PECO zone of the PJM power market [1]. PJM is a predominantly thermal system [2] and, as a result, spot prices vary widely along the day.

Mateus A. Cavaliere, Sergio Granville, Gerson C. Oliveira, Mario V.F. Pereira are with Power Systems Research, Rio de Janeiro, Brazil (e-mails: mateus@psr-inc.com, granville@psr-inc.com, gerson@psr-inc.com, mario@psr-inc.com, respectively).

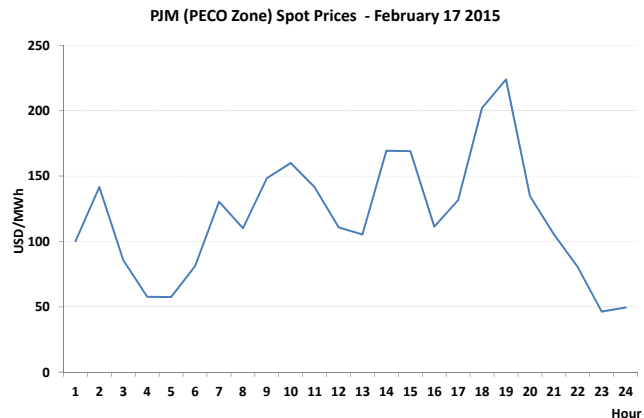

**Figure 1 - Hourly Spot Prices - PJM**

Figure 2 shows the average daily spot price, from 2007 to 2014, for the same system. Notice that, except for some spikes, prices show a quite stable behavior throughout the years.

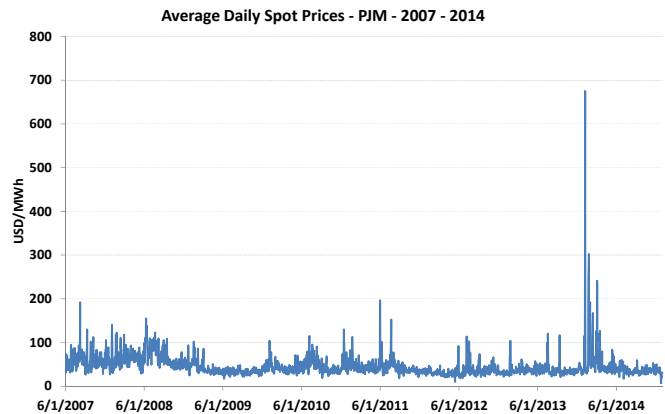

**Figure 2 - Average Daily Spot Prices**

Now consider a predominantly hydro system with large reservoirs, such as Brazil [3]. Spot prices in Brazil [4] are strongly correlated to the amount of water stored in the main hydro plant´s reservoirs, as shown in Figure 3. Notice how short-term price volatility is small, but medium to long-term volatility is huge.

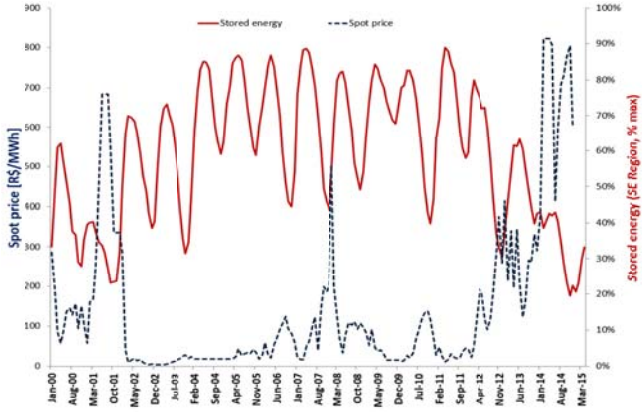

**Figure 3 - Brazilian Spot Prices x Stored Water**

Differences in spot prices across buses or zones are due to congestions in the transmission system, also presenting a great deal of volatility. Figure 4 shows the hourly spot price differences across three zones in the PJM market [1]. Notice the large gap between prices that occur during the day.

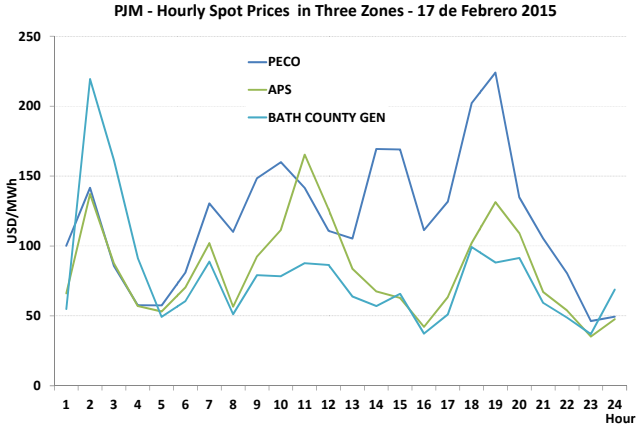

Figure 4 – Hourly Spot Price Differences Across Zones

One way to hedge against spot prices volatility is to establish future/forward contracts [4-7]. However, a great challenge for generator, load and trading companies in this framework is to define an optimal contract strategy in order to optimize their results, at the same time take into account their financial risk.

In this paper we present a methodology for future/forward electricity contract price projection which is a key element in the designing of optimal contract strategies.

The remainder of this work is organized as follows: in Section II we present an overview of the electricity contract marketing for generators, loads and trading companies together with the equilibrium concept. In Section III it is shown that competitive equilibrium in the contract market can be computed through welfare maximization. Section IV discusses extensions of the basic methodology. Section V presents applications of the methodology to the Brazilian electricity market and in Section VII the main conclusions are summarized.

## II. OVERVIEW OF THE ELECTRICITY CONTRACT MARKETING

A selling contract is a commitment to sell an amount of energy for a given time period and at given price by either generation or buying it in the spot market. Likewise, a buying contract is a commitment to buy an amount of energy for a given time period and at given price for consumption or selling in the spot market.

We suppose that agents have their supply and demand curve for contracts, associated to a given future time period (contract delivery period), as shown in Figure 5. For ease of presentation, we present in sections II and III a simple case of Equilibrium in Electricity Contract Market: one single time stage and just one generator and load, more complex cases will be discussed in the extensions.

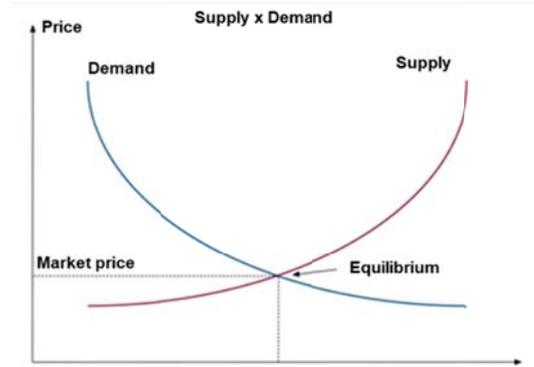

**Figure 5 - Supply and Demand Curves for Contracts**

Generator supply curve is built from the optimization of its risk adjusted expected net revenue for each contract price. The net revenue in the contract market is equal to

$$\hat{R}_g(q_s, p) = (\hat{g} - q_s) \times \hat{\pi} + q_s \times p \quad (1)$$

In equation (1) $\hat{g}$ is the generation amount, $q_s$ is the selling contract amount, $\hat{\pi}$ is the spot price and $p$ is the selling contract price[1]. The first component in equation (1) corresponds to the net revenue in the spot market and the second, the revenue in the contract market.

Likewise, load demand curve is built from the optimization of its risk adjusted expected net revenue for each contract price. The net revenue in the contract market is equal to

$$\hat{R}_d(q_b, p) = (q_b - \hat{d}) \times \hat{\pi} - q_b \times p \quad (2)$$

In equation (2), $\hat{d}$ is the load amount, $q_b$ is the buying contract amount and $p$ is the buying contract price. The first component in equation (2) corresponds to the net revenue in the spot market and the second, the cost in the contract market.

The risk adjusted expected value is defined as a convex combination of expected value and CVaR$_\alpha$ [8] for a given level of confidence level $\alpha$:

$$R_g^{Adj}(q_s, p) = \lambda_g E[\hat{R}_{gen}(q_s, p)] + (1 - \lambda_g)CVaR_\alpha[\hat{R}_{gen}(q_s, p)] \quad (3)$$

$$R_d^{Adj}(q_b, p) = \lambda_d E[\hat{R}_{load}(q_b, p)] + (1 - \lambda_d)CVaR_\alpha[\hat{R}_{load}(q_b, p)] \quad (4)$$

---

[1] In this work hats denote random variables

Weight $\lambda_g$, $\lambda_d$ is a risk aversion parameter for the generator and load.

For a given contract price $p$, generators sell contract amount associated to $p$, given by:

$$q_s(p) = argmax_{q_s}\{R_g^{Adj}(q_s, p)\} \quad (5)$$

Likewise, for given a contract price $p$ loads buy contract amount associated to $p$, given by

$$q_b(p) = argmax_{q_b}\{R_d^{Adj}(q_b, p)\} \quad (6)$$

Competitive equilibrium in the contract market corresponding to a pair $(p^0, q^0)$ of price and energy amount such that $q^0$ is the solution of both supply (4) and demand curve (5) associated to price $p^0$:

$$q^0 = q_s(p^0) = q_b(p^0) \quad (7)$$

In equation (7) $p^0$ is the future/forward equilibrium contract price, reflecting agent expectations about market conditions in the contract delivery date.

When uncertainties are represented by a discrete number of scenarios, using Rockafellar's representation of $CVaR$ [8], implies that problems (5)-(6) correspond to linear programming problems, parametrized by the price $p$. For a fix price $p$ the optimality conditions for problems (5) and (6) corresponds to a linear system of equality/inequality equations. However, when $p$ is considered as a variable the resulting system of equality/inequality equations become nonlinear due to products involving contract price and energy amount. This means that in principle computation of the competitive equilibrium price in the contract market requires solving a system of nonlinear equality/inequality equations. In the next section we will show that this equilibrium can be found through the resolution of a linear programing problem which corresponds to the welfare maximization of both agents.

## III. Equilibrium Price Computation in the Contract Market and Welfare Maximization

We start this section by observing that as the terms $q_s \times p$ and $q_b \times p$ in revenues equations for generators and loads are not random, expected value and $CVaR$ properties [9], imply that:

$$R_g^{Adj}(q_s, p) = \lambda_g E[(\hat{g} - q_s) \times \hat{\pi}] + (1 - \lambda_g)CVaR_\alpha[(\hat{g} - q_s) \times \hat{\pi}] + q_s \times p \quad (8)$$

$$R_d^{Adj}(q_b, p) = \lambda_d E[(q_b - \hat{d}) \times \hat{\pi}] + (1 - \lambda_d)CVaR_\alpha[(q_b - \hat{d}) \times \hat{\pi}] - q_b \times p \quad (9)$$

Note that in (8) $q_s \times p$ is equal to the revenue associated to contract selling. The other term corresponds to the negative of the costs in the spot market associated to contract energy delivering. Thus (8) corresponds to the seller's surplus. On the other hand, (9) corresponds to contract revenue in the spot market subtracted by the cost of purchasing it $(q_b \times p)$ which is equal to the buyer's surplus.

Supposing that uncertainties are represents by a discrete number of scenarios, Rockafellar's representation of CVaR [8] and (8) imply that problem (5) can be written as:

$$Max_{q_s}\left\{\lambda_g \frac{\sum_k R_{g,k}}{K} + (1 - \lambda_g)\left[a_g + \frac{\sum_k [R_{g,k} - a_g]^-}{K \times (1-\alpha)}\right] + q_s \times p\right\} \quad (10)$$

s.t.

$$R_{g,k} - (g_k - q_s) \times \pi_k = 0, k = 1, \ldots, K$$
$$q_s \geq 0$$

Where:

$K$    number of scenarios
$\pi_k$    Spot price in scenario $k$
$g_k$    Generation in scenario $k$
$R_{g,k}$    net revenue of generator in the spot market in scenario $k$
$[x]^-$    equals to $x$ if $x \leq 0$ and 0 if $x > 0$

Problem (8) is nonlinear. Performing a transformation, it can be rewritten as a linear programming problem:

$$Max_{q_s}\left\{\lambda_g \frac{\sum_k R_{g,k}}{K} + (1 - \lambda_g)\left[a_g + \frac{\sum_k y_{g,k}}{K \times (1-\alpha)}\right] + q_s \times p\right\} \quad (11)$$

s.t.      Dual variables:

$$R_{g,k} - (g_k - q_s) \times \pi_k = 0, k = 1, \ldots, K \quad \theta_{g,k}$$
$$q_s \geq 0 \quad \beta_g$$
$$y_{g,k} \leq 0, k = 1, \ldots, K \quad \gamma_{g,k}$$
$$y_{g,k} \leq R_{g,k} - a_g, k = 1, \ldots, K \quad \eta_{g,k}$$

Likewise, problem (6) with transformation above can be written as:

$$Max_{q_b}\left\{\lambda_d \frac{\sum_k R_{d,k}}{K} + (1 - \lambda_d)\left[a_d + \frac{\sum_k y_{d,k}}{K \times (1-\alpha)}\right] - q_b \times p\right\} \quad (12)$$

s.t.      Dual variables:

$$R_{d,k} - (q_b - d_k) \times \pi_k = 0, k = 1, \ldots, K \quad \theta_{d,k}$$
$$q_b \geq 0 \quad \beta_d$$
$$y_{d,k} \leq 0, k = 1, \ldots, K \quad \gamma_{d,k}$$
$$y_{d,k} \leq R_{d,k} - a_d, k = 1, \ldots, K \quad \eta_{d,k}$$

Where:

$d_k$    Load in scenario $k$
$R_{d,k}$    net revenue of load in the spot market in scenario $k$

Now suppose that there is a price $p^0$ which corresponds to equilibrium, that is:

$$q^0 = q_s(p^0) = q_b(p^0)$$

Optimality conditions for problems (11) at $q^0$, with price $p^0$, can be written as:

$$R_{g,k}^0 + q^0 \times \pi_k = \pi_k \times g_k, k = 1, \ldots, K \quad (13)$$
$$q^0 \geq 0 \quad (14)$$
$$y_{g,k}^0 \leq 0, k = 1, \ldots, K \quad (15)$$
$$y_{g,k}^0 - R_{g,k}^0 + a_g^0 \leq 0, k = 1, \ldots, K \quad (16)$$
$$\frac{\lambda_g}{K} - \theta_{g,k}^0 + \eta_{g,k}^0 = 0, k = 1, \ldots, K \quad (17)$$
$$(1 - \lambda_g) - \sum_k \eta_{g,k}^0 = 0 \quad (18)$$
$$\frac{(1-\lambda_g)}{K \times (1-\alpha)} - \gamma_{g,k}^0 - \eta_{g,k}^0 = 0, k = 1, \ldots, K \quad (19)$$
$$p^0 - \sum_k \pi_k \times \theta_{g,k}^0 - \beta_g^0 = 0 \quad (20)$$
$$\beta_g^0 \leq 0 \quad (21)$$
$$\gamma_{g,k}^0 \geq 0, k = 1, \ldots, K \quad (22)$$
$$\eta_{g,k}^0 \geq 0, k = 1, \ldots, K \quad (23)$$
$$\sum_k \theta_{g,k}^0 \times \pi_k \times g_k = \lambda_g \frac{\sum_k R_{g,k}^0}{K} + (1 - \lambda_g)\left[a_g^0 + \frac{\sum_k y_{g,k}^0}{K \times (1-\alpha)}\right] + q^0 \times p^0 \quad (24)$$

Equality (24) corresponds to the primal-dual equality condition.

Optimality conditions for problems (12) at $q^0$, with price $p^0$, can be written as:

$$R_{d,k}^0 - q^0 \times \pi_k = -\pi_k \times d_k, k = 1, \ldots, K \tag{25}$$
$$q^0 \geq 0 \tag{26}$$
$$y_{d,k}^0 \leq 0, k = 1, \ldots, K \tag{27}$$
$$y_{d,k}^0 - R_{d,k}^0 + a_d^0 \leq 0, k = 1, \ldots, K \tag{28}$$
$$\frac{\lambda_d}{K} - \theta_{d,k}^0 + \eta_{d,k}^0 = 0, k = 1, \ldots, K \tag{29}$$
$$(1 - \lambda_d) - \sum_k \eta_{d,k}^0 = 0 \tag{30}$$
$$\frac{(1-\lambda_d)}{K \times (1-\alpha)} - \gamma_{d,k}^0 - \eta_{d,k}^0 = 0, k = 1, \ldots, K \tag{31}$$
$$-p^0 + \sum_k \pi_k \times \theta_{d,k}^0 - \beta_d^0 = 0 \tag{32}$$
$$\beta_d^0 \leq 0 \tag{33}$$
$$\gamma_{d,k}^0 \geq 0, k = 1, \ldots, K \tag{34}$$
$$\eta_{d,k}^0 \geq 0, k = 1, \ldots, K \tag{35}$$
$$-\sum_k \theta_{d,k}^0 \times \pi_k \times d_k = \lambda_d \frac{\sum_k R_{d,k}^0}{K} + (1-\lambda_d) \left[ a_d^0 + \frac{\sum_k y_{d,k}^0}{K \times (1-\alpha)} \right] - q^0 \times p^0 \tag{36}$$

Equality (36) corresponds to the primal-dual equality condition.

Now consider the following linear optimization problem:

$$Max_{q_s, q_b} \left\{ \lambda_g \frac{\sum_k R_{g,k}}{K} + (1 - \lambda_g) \left[ a_g + \frac{\sum_k [R_{g,k} - a_g]^-}{K \times (1-\alpha)} \right] + \lambda_d \frac{\sum_k R_{d,k}}{K} + (1 - \lambda_d) \left[ a_d + \frac{\sum_k y_{d,k}}{K \times (1-\alpha)} \right] \right\} \tag{37}$$

s.t.                                                    Dual variables:
$$R_{g,k} - (g_k - q_s) \times \pi_k = 0, k = 1, \ldots, K \qquad \theta_{g,k}$$
$$q_s \geq 0 \qquad \beta_g$$
$$y_{g,k} \leq 0, k = 1, \ldots, K \qquad \gamma_{g,k}$$
$$y_{g,k} \leq R_{g,k} - a_g, k = 1, \ldots, K \qquad \eta_{g,k}$$
$$R_{d,k} - (q_b - d_k) \times \pi_k = 0, k = 1, \ldots, K \qquad \theta_{d,k}$$
$$q_b \geq 0 \qquad \beta_d$$
$$y_{d,k} \leq 0, k = 1, \ldots, K \qquad \gamma_{d,k}$$
$$y_{d,k} \leq R_{d,k} - a_d, k = 1, \ldots, K \qquad \eta_{d,k}$$
$$q_s - q_b = 0 \qquad \delta_{g,d}$$

Optimality conditions associated to problem (37) are:

$$R_{g,k} + q_s \times \pi_k = \pi_k \times g_k, k = 1, \ldots, K \tag{38}$$
$$q_s \geq 0 \tag{39}$$
$$y_{g,k} \leq 0, k = 1, \ldots, K \tag{40}$$
$$y_{g,k} - R_{g,k} + a_g \leq 0, k = 1, \ldots, K \tag{41}$$
$$R_{d,k} - q_b \times \pi_k = -\pi_k \times d_k, k = 1, \ldots, K \tag{42}$$
$$q_b \geq 0 \tag{43}$$
$$y_{d,k} \leq 0, k = 1, \ldots, K \tag{44}$$
$$y_{d,k} - R_{d,k} + a_d \leq 0, k = 1, \ldots, K \tag{45}$$
$$q_s - q_b = 0 \tag{46}$$
$$\frac{\lambda_g}{K} - \theta_{g,k} + \eta_{g,k} = 0, k = 1, \ldots, K \tag{47}$$
$$(1 - \lambda_g) - \sum_k \eta_{g,k} = 0 \tag{48}$$
$$\frac{(1-\lambda_g)}{K \times (1-\alpha)} - \gamma_{g,k} - \eta_{g,k} = 0, k = 1, \ldots, K \tag{49}$$
$$-\delta_{g,d} - \sum_k \pi_k \times \theta_{g,k} - \beta_g = 0 \tag{50}$$
$$\beta_g \leq 0 \tag{51}$$
$$\gamma_{g,k} \geq 0, k = 1, \ldots, K \tag{52}$$
$$\eta_{g,k} \geq 0, k = 1, \ldots, K \tag{53}$$
$$\frac{\lambda_d}{K} - \theta_{d,k} + \eta_{d,k} = 0, k = 1, \ldots, K \tag{54}$$
$$(1 - \lambda_d) - \sum_k \eta_{d,k} = 0 \tag{55}$$
$$\frac{(1-\lambda_d)}{K \times (1-\alpha)} - \gamma_{d,k} - \eta_{d,k} = 0, k = 1, \ldots, K \tag{56}$$
$$\delta_{g,d} + \sum_k \pi_k \times \theta_{d,k} - \beta_d = 0 \tag{57}$$
$$\beta_g \leq 0 \tag{58}$$
$$\gamma_{g,k} \geq 0, k = 1, \ldots, K \tag{59}$$
$$\eta_{g,k} \geq 0, k = 1, \ldots, K \tag{60}$$
$$\sum_k \theta_{g,k} \times \pi_k \times g_k - \sum_k \theta_{d,k} \times \pi_k \times d_k = \lambda_g \frac{\sum_k R_{g,k}}{K} + (1 - \lambda_g) \left[ a_g + \frac{\sum_k y_{g,k}}{K \times (1-\alpha)} \right] + \lambda_d \frac{\sum_k R_{d,k}}{K} + (1 - \lambda_d) \left[ a_d + \frac{\sum_k y_{d,k}}{K \times (1-\alpha)} \right] \tag{61}$$

Comparing the set of linear equality/inequality (13)-(36) with the set (38)-(61) we conclude that the optimal solutions of problems (11) e (12) corresponding to an equilibrium, satisfies optimality conditions of problem (37) with the following correspondence

$$q_s = q_b = q^0 \tag{62}$$
$$\delta_{g,d} = -p^0 \tag{63}$$

This implies that the solution of an equilibrium problem corresponds to the solution of the linear programming problem (37) where the dual of the constraint that specify that the energy amount of selling contract is equal to the energy amount of selling contract, is equal to the equilibrium price. Conversely a solution of the linear programming problem (37) corresponds to an equilibrium in the contract market with the equilibrium price equal to the dual of the constraints that specify that the energy amount of selling contract is equal to the energy amount of selling contract.

Now note problem (37) corresponds to the maximization of the sum of the objective functions associated to the contract supply and demand curves (terms $q_s \times p$ and $q_b \times p$ cancel out due to constraint $q_s - q_b = 0$ of the problem). As these objective functions are equal to the producer surplus and consumer surplus, problem (37) corresponds to the economic surplus or total welfare maximization (see Figure 6)

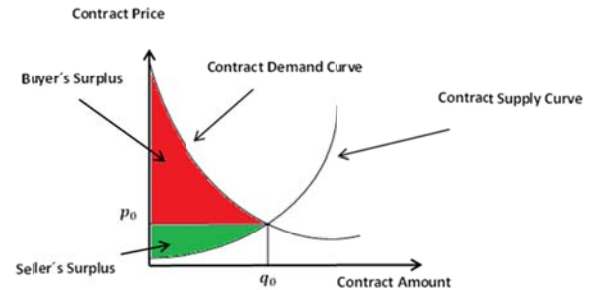

**Figure 6 – Welfare optimization**

## IV. EXTENSIONS OF THE METHODOLOGY

In this section we will consider several extensions of the basic methodology.

### A. More than two agents

The development in preceding sections can be extended to consider equilibrium in the contract market for a set of generators, loads and trading companies. Pure trading companies do not have any physical assets (generators or loads) but can establish selling or buying contracts. In this case its net revenue is equal to:

$$\hat{R}_t(q_s, q_b, p) = (q_b - q_s) \times \hat{\pi} + (q_s - q_b) \times p \qquad (64)$$

In equation (64) $q_s$, $q_b$, represent selling and buying contract amounts; the first component corresponds to the net revenue in the spot market and the second, the net revenue in the contract market.

Based on the optimization of its risk adjusted expected revenue, each agent establishes in this case its selling/buying contract amount as a function of contact price. Equilibrium in the contract market corresponds now to a contract price $p^0$, associated to energy amount for a set to buying and selling contracts $\{q_{b,i}^0, i = 1, \ldots, I_b\}$, $\{q_{s,i}^0, j = 1, \ldots, J_s\}$ contracts such that:

$$\sum_i q_{b,i}^0 = \sum_j q_{s,i}^0 \qquad (65)$$

The extension of the methodology developed in previous sections to accomplish this case is immediate.

### B. More than a time period

Electricity contract may involve a few hours for a given month (e.g. peak hours for a calendar month contracts in Nymex), days (peak days for a calendar month contracts in Nymex), weeks (e.g. weekly contracts in Nordpool) months, quarters, years (e.g. monthly, yearly contracts in Nordpool and Brazil), etc.

In this case contract energy amount is indexed by time period ($m$) (e.g. $q_{s,m}$ or $q_{b,m}$); it has an associated shape, for instance:

$$q_{s,m} = q_s \times v_m, \ \mathrm{m} = 1, \ldots, \mathrm{M} \qquad (66)$$

where

$q_s$   Decision variable – contract reference energy amount

$v_m$   Vector associated to contract energy distribution along its horizon

Also, spot prices, generation and loads are indexed by time period:

$$\{\pi_{m,k}, g_{m,k}, d_{m,k}, m = 1, \ldots, M; \ k = 1, \ldots, K\} \qquad (67)$$

As a result, agent net revenues and market settlements are established for each time period.

The extension of the methodology developed in previous sections to accomplish this case is immediate.

### C. Time Variations and Uncertainties in forward contract prices

As mentioned before, forward prices reflect agent expectations about future market condition; these expectations change over time. In order to model these features we have to consider the dynamic behavior of the key random variables that affects agent revenues. Usually the distributions of these variables in successive time periods are not independent. Their evolution can be represented by a tree structure which reflects transitions between system states along successive time periods (see Figure 7).

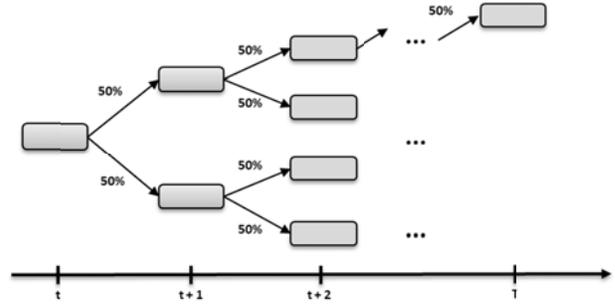

**Figure 7 - Tree structure**

There are many ways to generate a tree and its associated transition probabilities from an underlining multi-stage stochastic process. One of them is through an approximation of the stochastic process distribution probability [10]. However, as in here we do not have a parametric distribution associated to the stochastic process, it is adopted an alternative approach based on a sample estimation procedure as described next. First of all, as the key random variables associated to energy markets correspond to time series, their representation through a finite set of scenarios (sample) corresponding to a finite number of trajectories. This sample of trajectories may be generated through a hydrothermal optimization model or a stochastic model adjusted to historical data.

The first step in the procedure is to define the tree topology in terms of a set of nodes for each time stage and direct edges. Each edge has a From node and a To node, and connects a node of one time stage (From node) to a node in the next time stage (To node). A root is a node which does not corresponds to a To node of any edge. Each node (except for the root) has just one preceding node, that is the number of nodes in each time stage is greater or equal the number of nodes in the preceding time stage.

Once the tree topology is defined, a variable or set of variables is selected for the clustering technique. The next step is to apply the clustering technique and compute the transition probabilities as described next.

For each time stage the values associated to the set trajectories in that time stage is clustered into $nb_t$ subsets, where $nb_t$ is the number of tree nodes in that time stage.

The transition probabilities are computed as follows: let $k$ a node in time stage $t$ and $j$ a node in time stage $t + 1$, the sample estimation for transition probability between these two nodes is equal to $nsb_{t,t+1,j}/nsb_{t,k}$, where $nsb_{t,k}$ is the number of trajectories whose values belong to cluster $k$ at stage $t$ and $nsb_{t,t+1,j}$ is the number of trajectories whose values belong to cluster $k$ in stage $t$ and to cluster $j$ in stage $t + 1$.

Given the construction above, note that besides tree transition probabilities we have, within each node, scenarios which represent uncertainties associated to the key random variables that affects agent revenues. As a result, the methodology developed above can be extended to compute forward contract price associated to future delivery period and system state related to each node. This feature allows us to compute for a forward contract, associated to a future delivery period, the evolution of its price probability distribution overtime; this is important to quantify financial risks associated to forward contracts and also the Value of Waiting.

## V. Case Examples

For the case examples illustrated in this paper 1200 spot price scenarios are considered along the years of 2022 to 2026 for the Brazilian system. Figure 8 shows aggregate spot price distribution for the year of 2022, 2022-2024 and 2022-2026. It can be seen that there are many scenarios with low spot prices and a small number of scenarios with high spot prices, that is, spot price distribution is highly asymmetrical.

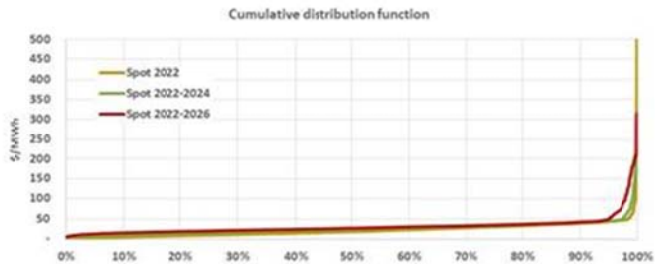

**Figure 8 - Spot Price Distribution**

Three types of simulation will be discussed. In the first the effect of generator capacity is analyzed; the second concerns the effect of risk aversion parameter ($\lambda$) and the third is related to uncertainties in future contract price.

All associated optimization problems are solved by the Xpress software [11-12].

### A. Effect of generation capacity

In this case the risk aversion parameter for both generator and load was set to 0.8, and three types of contracts were considered: the first one with one year duration and delivery period year of 2022, the second with three years duration and delivery period from the year of 2022 to 2024 and the third with five years duration and delivery period from the year of 2022 to 2026. Load size is 100 MW and generator capacity varies from 100 MW to 130 MW (overcapacity ranging from 0% to 30%).

Figure 1Figure 9 shows contract equilibrium price corresponding to generator overcapacity ranging from 0% to 30% and average spot price for the years of contract delivery period.

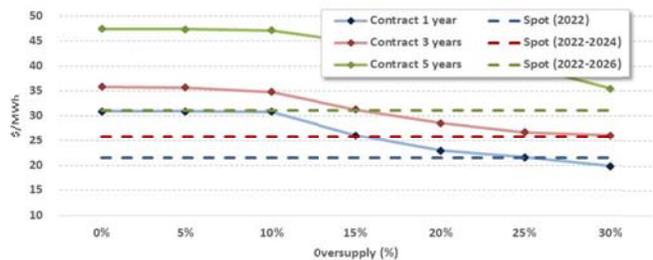

**Figure 9 - Effect of Energy Overcapacity**

First, note that contract price tends to be higher than average spot price along the corresponding delivery period. This is mainly due to high spot scenarios together with its distribution asymmetry (see Figure 8) which imply in a risk premium charged by the contract seller.

In addition, starting in 2022, contract price increases with contract durations. This is mainly due spot prices increasing behavior from 2022 to 2026.

Besides, it can be observed that equilibrium contract price decreases as a function of generator overcapacity. This happens because, as overcapacity becomes higher, an increasing excess of energy is supplied to the contract market and this implies that the generator tends to accept a lower contract price.

### B. Effect of risk aversion parameter

Now let us examine the impact of the risk aversion parameter on contract price. In this case it was considered one year contract for the year of 2022, a load of 100 MW, generator capacity of 110 MW, three values for the load risk aversion parameter ($\lambda_q = 1.0, 0.5, 0.0$) and three values for the generator risk aversion parameter ($\lambda_g = 1.0, 0.5, 0.0$).

Figure 10 shows the resulting contract equilibrium price. First, note that when the load or generator is risk neutral ($\lambda_q = 1$ or $\lambda_g = 1.0$) the contract price is equal to average spot price for the year of 2022 ($\approx$ \$22). This happens because, if load is risk neutral, it would does not have any motivation to accept a risk premium charged by the generator due spot asymmetric distribution. By the same reasoning, if generator is risk neutral it does not have any motivation to charge a risk premium from the load.

As the load has some level of risk aversion ($\lambda_q < 1$) contract price increases with generator level of risk aversion. This happens because generation charged risk premium increases with its level of risk aversion and this risk premium is accepted by the load, due to its risk aversion. Also, as the generator has some level of risk aversion ($\lambda_g < 1$) the contract price increases with load level risk aversion. As load level of risk aversion increases it will be increasingly willing to accept generator risk premium.

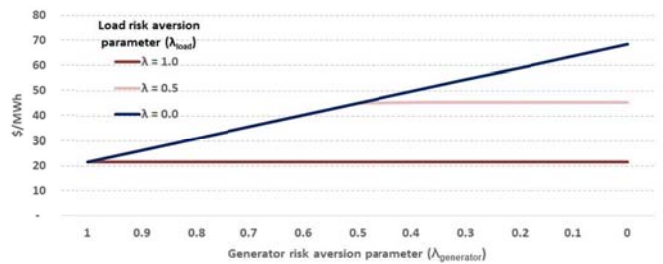

**Figure 10 – Effect of Risk Aversion Profile**

### C. Uncertainties in future contract prices

Now we turn to the computation of the future contract price probability distribution. Consider a one-year duration contract with delivery year of 2025 and that we want to study de evolution of its price from the beginning of 2022 to the beginning of 2025. For this purpose, consider a tree with the topology shown in Figure 11. At the beginning there were 1200 spot price scenarios which correspond to trajectories in 2022. Then based on spot prices in 2023 these trajectories were clustered into two classes containing 600 values each (see Figure 11). Afterwards, based on spot prices in 2024, each of the 600 trajectories were clustered into two classes in 2024, a total of four classes at that year. Finally, based on spot prices in 2025, each of the 300 trajectories were clustered into two classes in 2025, a total of eight classes in 2025.

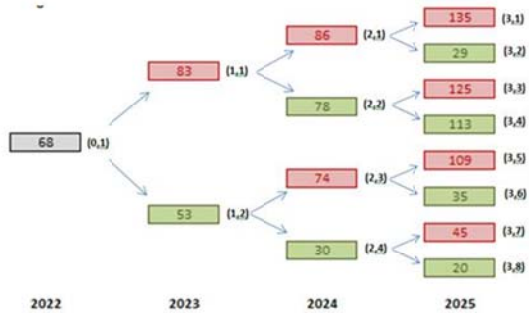

**Figure 11 – Tree Topology and Results**

Considering variables that affect agent revenues in 2025 associated to each of the 1200 trajectories, the methodology can be applied to compute forward price contract for that year. This is the price ($ 68) in 2022 for the forward contract delivery in 2025 as shown in the figure. Next, applying the methodology to the variables that affect agent revenues in 2025 which belong to the 600 trajectories associated to the tree node (1,1), future contract price for 2025 in that state is $ 83 as shown in the figure. Doing the same for the other 600 trajectories associated to the tree node (1,2), future contract price for 2025 in that state is $ 53 as shown in the figure. Applying the same procedure for the other nodes future contract price for 2025 are obtained for each of the remaining states as shown in the figure.

As the 1200 trajectories are equally likely, contract price for 2025 in 2022 is $68 but there is a 50% probability that it will be $83 in 2023 and 50% that it will be $53. For 2024, the future contract price for 2025 could be $86, $78, $74 and $30 with a 25% probability associated to each value. Likewise, there are eight possible values for future contract price at the beginning of 2025 to deliver energy in 2025, each associated with a 12.5% probability.

From the above it can be seen that contract prices for a future delivery period not only changes overtime but its probability distribution also changes.

Uncertainties in contract prices imply in risks for buying and selling agents in the contract market. For instance, consider in the example, a selling contract, established in 2022 at a price $68 for delivery energy in 2025. Suppose that market state in 2024 is 1 (node (2,1) in the tree). Contract price for delivery energy in 2025 in that state is $86. The selling agent who established the contract in 2022 had a loss because in that state contract price for energy delivery in 2025 is higher. That is, the value of the contract established in 2022 is in that state -$18 ($68 - $86). Considering all states, the value of a selling contract established in 2022 for energy delivery in 2025, can be in 2024 equal to -$18, -$10, -$6 or $38, with a 25% probability associated to each value.

## VI. CONCLUSIONS

In this paper it was presented a methodology to estimate forward contract price in energy market. The methodology is based on the computation of equilibrium in a competitive market, where market agents (buyers and sellers) submit their corresponding supply and demand curves for contracts.

Contract supply and demand curves are based on the optimization of each agent risk adjusted expected revenue established in terms of a convex combination of the expected revenue and CVaR.

When uncertainties are represented by a discrete number of scenarios, a key result in the paper is that contract price corresponds to the dual variable of the equilibrium constraints in the linear programming problem associated to the optimization of total agent welfare.

Besides computing equilibrium contract price for a given year, the methodology can also be used to compute the evolution of the price distribution associated to a contract with a future delivery period; this is an important topic to quantify forward contract risks.

Numerical examples provided in the paper show the attractiveness of the methodology to deal with forward contract pricing, a key element in the designing of optimal contract strategies.

## VII. ACKNOWLEDGMENT

The authors gratefully acknowledge the contributions of Rodrigo Cavalcanti and Bernardo Bezerra from PSR on the development of this work.

## IX. BIOGRAPHIES

**Mateus A. Cavaliere**, has a BSc in Electrical Engineering from UFRJ. Between 2014 and 2015 he, studied in Energy Systems at Hochschule Ruhr

West, Germany. Joined PSR in 2015 and has been working in (i) tariffs studies and energy pricing in the environments of regulated and free contracting; (ii) strategy assessment studies for supplying energy to consumers; (iii) economic and financial evaluations of generation, distribution and transmission companies; and (iv) regulatory advisory.

**Sérgio Granville** has a BSc and a MSc degrees in Applied Mathematics from PUC/RIO and a PhD degree in Optimization from Stanford University. He has coordinated software development and optimization studies for portfolios involving physical and financial assets, financial evaluation of generation projects, strategic bidding in energy auctions, pricing of energy contracts and options, and the development of mathematical tools for the analysis of regulatory aspects. Dr. Granville also participated in the investment analysis of generation/transmission projects and contracting strategies, both in Brazil and in foreign countries Before joining PSR, he was a project manager at Cepel, where he coordinated research and developments of optimization tools in the areas of optimal power flow and VAr planning. Dr. Granville has authored and co-authored more than fifty papers in refereed journals and conference proceedings.

**Gerson C. Oliveira** has a BSc degree in EE from PUC/Rio and PhD degree in Optimization from COPPE/UFRJ. He has coordinated software development for transmission planning software, stochastic streamflow modeling, supply reliability evaluation. He authored more than fifty papers in refereed journals and conference proceedings, as well as chapters for books on power systems planning and operations. He is a member of the Brazilian Automatic Control Society. Previously he worked at Cepel, where he was a project manager in the areas of optimization and statistical techniques applied to power systems planning and operations.

**Mario Veiga Pereira** has a BSc degree in EE from PUC/Rio and PhD degree in Optimization from COPPE/UFRJ. He founded PSR in 1987 and has been the company's CEO until December 2018, when he became Chief Innovation Officer (CINO). He was the lead formulator of the country's energy contracting auctions (80 GW of new generation capacity contracted since 2005, for about US$ 550 billion); co-designed the country's bioelectricity program (sugarcane biomass); and developed methodologies for wind energy auctions. Some of his recent international consulting activities include the development of a new generation planning model for the US Pacific Northwest and Western Interconnection; a price forecasting system for the Scandinavian Nordpool; the power market reform of Turkey, Vietnam and Mexico; the interconnection study of 8 South American countries; and studies for renewable integration in several countries worldwide. Dr. Pereira is a leading developer of advanced optimization methods. His stochastic dual dynamic programming (SDDP) algorithm is a worldwide reference and applied in dozens of countries. He also developed novel methods for optimal expansion planning and supply reliability evaluation. Dr. Pereira is an IEEE Fellow; an elected member of Brazil's Academy of Sciences and National Academy of Engineering; was awarded a Presidential Medal for his contribution to Brazil's electricity sector; a Scientific Merit Medal for his research contributions; and is a co-recipient of the Franz Edelman Award for the development of stochastic optimization tools for the operation of hydrothermal systems. He is the author and co-author of four books and about 200 papers in refereed journals